\magnification 1200

\documentstyle{amsppt}
\loadbold
\NoBlackBoxes
\def\dbend{{\manual\char127}} 
\def\d@nger{\medbreak\begingroup\clubpenalty=10000
  \def\par{\endgraf\endgroup\medbreak} \noindent\hang\hangafter=-2
  \hbox to0pt{\hskip-\hangindent\dbend\hfill}\ninepoint}
\outer\def\danger{\d@nger}
\def\dd@nger{\medbreak\begingroup\clubpenalty=10000
  \def\par{\endgraf\endgroup\medbreak} \noindent\hang\hangafter=-2
  \hbox to0pt{\hskip-\hangindent\dbend\kern1pt\dbend\hfill}\ninepoint}
\outer\def\ddanger{\dd@nger}
\def\today{\ifcase\month\or January\or February\or March\or
April\or May\or June\or July\or August\or September\or October\or
November\or December\fi \space\number\day, \number\year}

\def\qed{\hfill$\square$}

\def\ju{\vskip .8truecm plus .1truecm minus .1truecm}

\def\sju{\vskip .4truecm plus .1truecm minus .1truecm}
\def\lra{\longrightarrow}

\def\qq{\lq\lq}
\def\lqq{\qq}
\def\pr{^\prime}

\def\ets{\emptyset}
\def\se{\subseteq}


\def\shimply{\Rightarrow} 

\def\NN{\Bbb N}

\def\al{\alpha}
\def\ep{\epsilon}

\def\be{\beta}
\def\ga{\gamma}

\def\sst#1#2{\{ {#1}:{#2} \}}
\def\Ker{\hbox{\rm Ker\,}}

\def\ima{\hbox{\rm Im\,}}

\def\invlim{\varprojlim}

\def\noi{\noindent}




\hsize 12truecm
\vsize 20.8truecm

\nologo


\parindent=18pt
\document
\noindent{\bf International Journal of Pure and Applied Mathematics}
\hrule
\noi{\bf Volume 10\quad No.3\quad 2004, 349--356}
\vskip 2truecm
\ju
\centerline{\bf A NOTE ON SURJECTIVE INVERSE SYSTEMS}
\footnote[]{\noi Received December 1, 2002\qquad \hfill \copyright 2004 Academic Publications}
\ju
\centerline{Radoslav M. Dimitric}
\centerline{ Texas A\&M University, Department of Mathematics}
\centerline{ PO Box 1675, Galveston, TX 77553, USA}
\centerline{e-mail: dimitric\@tamug.edu}
\leftheadtext{R.M. Dimitri\'c}
\rightheadtext{A note on surjective inverse systems}
\vskip 1.5truecm
\noi {\bf Abstract:}
Given an upward directed set $I$ we consider surjective
$I$-inverse systems $\{X_\al,f_{\al\be}:X_\be\lra X_\al|
\al\leq\be\in I\}$,
namely those inverse systems that have all $f_{\al\be}$
surjective. A number of properties of $I$-inverse systems have
been investigated; such are the Mittag-Leffler condition,
investigated by Grothendieck and flabby and semi-flabby 
$I$-inverse systems studied by Jensen.  We note that flabby implies
semi-flabby implies surjective implies Mittag-Leffler. Some of the
results about surjective inverse systems have been known for some
time. The aim of this note is to give a series of equivalent 
statements and implications involving surjective inverse systems 
and the systems satisfying the Mittag-Leffler condition,
together with improvements of established results, as well as
 their relationships with the already known, but scattered facts. 
The most prominent results relate cardinalities of the index sets 
with right exactness of the inverse limit functor and the non-vanishing
of the inverse limit -- connections related to cohomological dimensions.

\noi{\bf AMS Subject Classification:} 18A30, 18G10, 18G20, 16E10, 13D05

\noi{\bf Keywords:} Surjective inverse system, Mittag-Leffler condition, $G$-sets,
surjective cohomological dimension

\sju
\pageno=349
For an (abelian) category $\Cal C$ with infinite products, and an 
(upward directed) ordered index set $I$,
an $I$-inverse system $\{ X_i, f_{ij}\}$ in $\Cal C$ is said to
satisfy the {\it Mittag-Leffler condition} (or the {\it ML
condition}, for short),
if $\forall i\in I$ $\exists j\geq i$ such
that, $\forall k\geq j$ $f_{ik}X_k=f_{ij}X_j$, in other words, if
for every $i\in I$, the decreasing sequence $\{ f_{ik}X_k\}_{k\in I}$
of submodules of $X_i$ stabilizes from certain index on; we will
also say that $\{ X_i, f_{ij}\}$ is {\it eventually stable}.
Note that for any $i\in I$, $\{ f_{ij}X_j\}_{j\geq i}$ is a
 decreasing family  of subobjects of $X_i$. Then the {\it
object of universal images} $X\pr_i=\inf_{j\geq i}f_{ij}X_j
(=\cap_{j\geq i}f_{ij}X_j)$ of
$X_i$ exists and, for $f_i:X=\invlim X_i\lra X_i$ -- the
canonical morphism,  $f_iX\se X\pr_i$ with $f_{ij}X\pr_j\se
X\pr_i$, for $j\geq i$. This makes $\{ X\pr_i, f_{ij}|X\pr_j=f\pr_{ij}\}$
an $I$-inverse system (the {\it  inverse system of universal images}), 
with $\invlim X\pr_i=X$. If ML holds, for
$\{ X_i, f_{ij}\}$ 
then the restriction $f\pr_{ij}=f_{ij}|X\pr_j$ is surjective, for 
all $j\geq i$
(there is a $k\geq j$ such that, for every $l\geq k\geq j$, 
$f_{jl}X_l=f_{jk}X_k$, hence
$X_j\pr=f_{jk}X_k$ and $f_{il}X_l=f_{ik}X_j$, thus we have
$X_i\pr=f_{ik}X_k=f_{ij}X_j\pr$). 

We now  construct  a  useful example of a surjective $I$-inverse
system  $\sst{E_\al,\epsilon_{\al\be}}{\al\leq\be\in I}$ of non-empty sets, 
for every non-empty 
upward directed set $I$. The case when $I$ has a maximal (the
maximum) element is usually favorable in considerations about
$\invlim$, hence we treat a more difficult case when $I$ has no 
maximal elements.

For $\al\in I$ let $E_\al$ denote the set of ordered
even-tuplets $(\al_1, \al_2,\dots,\al_{2n-1},\al_{2n})$ of elements of
$I$, with  the following properties: 
\roster
\item $\al_{2n-1}=\al$,
\item $\forall\,\, 1\leq i\leq n$, \,\,  $\al_{2i-1}\leq\al_{2i}$, and 
\item $\forall\,\,  1\leq j<i\leq n$, \,\, $\al_{2i-1}\not\leq\al_{2j-1}$. 
\endroster
All $E_\al$ are non-empty,
by non-maximality; they are also disjoint. To construct the maps, let $\al\leq\be$ and
$(\be_1,\dots,\be_{2m-1},\be_{2m})\in E_\be$, let $j$ be the smallest
numeral such that $\al\leq\be_{2j-1}$; define
$\ep_{\al\be}:E_\be\lra E_\al$ by
$\ep_{\al\be}(\be_1,\dots,\be_{2m-1},\be_{2m})=
(\be_1,\dots,\be_{2j-2},\al,\be_{2j})$. 
To show that $\{ E_\al, \ep_{\al\be}\}$ is a surjective
inverse system, start with an
$x=(\al_1,\dots,\al_{2n})\in E_\al$ 
(hence $\al_{2n-1}=\al$). Pick any $\ga>\be\geq \al$; we prove that 
$\ep_{\al\be}(y)=x$,
where $y=(\al_1,\dots,\al_{2n},\be,\ga)$. The non-trivial point to prove
is that $y\in E_\be$: $x\in E_\al$ satisfies conditions (1)--(3),
thus we need only prove that $\be\not\leq \al_{2l+1}$, for
$l=1,2,\dots, n-1$. This must be the case, for if $\be\leq
\al_{2l+1}$, then $\al<\be$ would imply $\al_{2n-1}=\al<\al_{2l+1}$, but this
is prevented by condition (3) for $x$. The verification that
$\ep_{\al\be}(y)=x$ is straightforward.

$G$-{\it sets}:
Given a  group $G$, a non-empty set $X$ is called a
$G$-{\it set}, if $G$ operates on $X$, i.e. if there is an operation 
$G\times X\lra X$, $(g,x)\mapsto gx$, such that 
$(g_1g_2)x=g_1(g_2x)$, and $ex=x$; a
consequence is that if $gx=y$, then $x=g^{-1}y$. A non-empty
$G$-set is {\it transitive}, if $\forall x_1, x_2\in X$, there
exists a $g\in G$ with $gx_1=x_2$; if another element $g_1$
satisfies the same equation (i.e. $g_1x_1=gx_1$), then $g^{-1}g_1$ 
belongs to the
{\it isotropy subgroup (stabilizer)} $H=\sst{g\in G}{gx_1=x_1}$ of $x_1$. 
If $G$ is the trivial group, then every set is trivially a
$G$-set, but non-transitive in general. Every group $G$ is a
transitive $G$-set.

\proclaim{\quad Theorem 1}  
For a non-empty upward directed set $I$, the following 
are equivalent:
\roster
\item $I$ has a maximal element, or it contains a countable cofinal
sequence. 

\item Every surjective $I$-inverse system
$\sst{X_\al}{f_{\al\be}}$ of non-empty sets has a
non-empty inverse limit. 

\item For every surjective map  
$\bold g=(g_\al)_{\al\in I}:\sst{E_\al}{\ep_{\al\be}}\lra
\sst{S_\al}{\sigma_{\al\be}}$, 
of $I$-inverse systems of sets, such that all
$\ep_{\al\be}$ are surjective and $\sigma_{\al\be}$ are injective, 
the induced inverse limit map $\invlim\bold g:\invlim{\bold E}\lra
\invlim{\bold S}$ is likewise surjective.

\item Every $I$-inverse system of non-empty sets $\sst{X_\al}{f_{\al\be}}$
that satisfies the ML condition has a non-empty inverse limit. 

\item  For every group $G$, every surjective $I$-inverse system 
of non-empty transitive 
$G$-sets has a non-empty inverse limit.
\endroster
\endproclaim

{\sl Proof.} (1)$\shimply $(2): If there is a maximum $\al_0\in
I$, then pick an $x_{\al_0}\in X_{\al_0}$ and define
$x_\al=f_{\al\al_0}(x_{\al_0})$, for every $\al\leq\al_0$; then
$(x_\al)_{\al\in I}\in\invlim X_\al\neq\ets$ (notice that
surjectivity of the inverse system is not needed here). 
If $J=\sst{\al_n}{n\in\NN}$ is a 
countable cofinal
subset of $I$, then $\invlim_IX_\al\cong\invlim_JX_\al$.
Pick an $x_{\al_1}\in X_{\al_1}$; by surjectivity there is an
$x_{\al_2}\in X_{\al_2}$ with
$f_{\al_1\al_2}(x_{\al_2})=x_{\al_1}$. Thus we construct
inductively $(x_\al)_{\al\in J}\in\invlim_J X_\al\neq\ets$.

(2)$\shimply$(3): Let $(s_\al)_{\al\in I}\in\invlim S_\al$;
this means that
$\forall \al\leq\be, \sigma_{\al\be}(s_\be)=s_\al$. Define 
$E\pr_\al=g^{-1}_\al(s_\al)\neq\ets$. We have
$\ep_{\al\be}(E\pr_\be)\se E\pr_\al$. This is because
$g_\al(\ep_{\al\be}g^{-1}_\be(s_\be))=\sigma_{\al\be}g_\be
g^{-1}_\be s_\be=\sigma_{\al\be}s_\be=s_\al$. Hence, for
$\ep\pr_{\al\be}=\ep_{\al\be}|E\pr_\al$, we have an $I$-inverse
system  $\{ E\pr_\al,\ep\pr_{\al\be}\}$ of non-empty sets. Moreover,
it is a surjective system, for if $y\in E\pr_\al$, surjectivity
of $\ep_{\al\be}$ ensures existence of an $x\in E\pr_\be$ with
$\ep_{\al\be}(x)=y$. This $x$ is in $E\pr_\be$, since
$s_\al=\sigma_{\al\be}s_\be$ and
$s_\al=g_\al\ep_{\al\be}(x)=\sigma_{\al\be}g_\be(x)$, hence
$\sigma_{\al\be}s_\be=\sigma_{\al\be}g_\be(x)$ and injectivity of
$\sigma_{\al\be}$ ensures $s_\be=g_\be(x)$. 

(3)$\shimply$(1): If $I$ has a maximum, there is nothing to
prove. Otherwise, consider a map $\bold g:\bold E\lra\bold S$ between
the special $I$-inverse systems $\bold E=\{E_\al,f_{\al\be}\}$
constructed in the introduction and
$\bold S=\{S_\al,\sigma_{\al\be}\}$, $S_\al=\{\al\}$,
$\sigma_{\al\be}(\be)=\al$. Both of the systems
are surjective and $\sigma_{\al\be}$ are bijections.
Define $g_\al(\al_1,\al_2,\dots,\al,\al_{2n})=\al$;  
this $\bold g$ is clearly surjective, and by the assumption, $\invlim
\bold g:\invlim \bold E\lra \invlim \bold S=(\al)_{\al\in I}$ 
is also surjective,
thus there is an $(e_\al)_{\al\in I}\in\invlim E$ that maps to
$(\al)_{\al\in I}$. These $e_\al$'s will produce a desired,
cofinal sequence in $I$ as follows:
Looking into the set of ending coordinates of all the $e_\al$'s, we
see that that set is cofinal in $I$ since, for every $\al\in I$,
$\al=\al_{2n-1}<\al_{2n}$. This means that we would prove the claim 
if we show that this set either has a maximal
element, or forms a countable sequence. We are assuming that
there is no maximal element. Note that if
$e_\al$ and $e_\be$ are tuplets of the same length $2n$, then $\al=\be$,
since if $\ga$ is chosen so that $\ga>\al,\be$, then $f_{\al\ga}e_\ga=e_\al$
and $f_{\be\ga}e_\ga=e_\be$; then the definition of the inverse system
morphisms $f_{ij}$
implies that the ending coordinates of $e_\al$ and $e_\be$ are
certain coordinates $\ga_{2l}$ and $\ga_{2m}$ of $e_\ga$. 
By the assumption of same length, $l=n=m$ and the ending coordinate 
of $e_\al$ is
$\ga_{2m}=\ga_{2l}=$ the ending coordinate of $e_\be$. The sizes
of all the tuplets $e_\al$ are not bounded, for otherwise their
ending coordinates would form a finite cofinal subset of $I$, 
hence it would have a maximal
element; since the lengths are not bounded, consider the
countable sequence of ending coordinates of each even-tuplet as
the desired cofinal sequence.  

(4)$\shimply$(2) holds since surjective inverse systems are special
cases of the ML systems.

For (2)$\shimply (5)$: Use the obvious forgetful functor.

(2)$\shimply$(4):
We have
already mentioned, that the $I$-inverse system of universal images
$f\pr_{\al\be}:X\pr_\be\lra X\pr_\al$ 
(where $X\pr_\be=f_{\be,\be+1}(X_{\be+1})$ -- we can assume that $I=\NN$, without
loss of generality)
is a surjective inverse system with $\invlim X\pr_\al=\invlim X_\al$,
provided $\{X_\al, f_{\al\be}\}$ satisfies the ML condition;
since $X_\al$'s are non-empty, then also $X\pr_\al$'s are
non-empty. By (2) $\invlim X\pr_\al\neq\ets$ and the claim is established.

(5)$\shimply$(1): 
Let $G$ be the (additively written) free abelian group on a set 
of generators $g_{ij}$
($i\leq j$ in $I$), and for each  $\al\in I$ let  $H_\al$ be the subgroup 
of  $G$  generated
by the elements 

$(a)\hfill    g_{ij} + g_{jk} - g_{ik},\quad k > j > i\geq\al. \hfill$ 

Let us define  $X_\al$  to be a transitive $G$-set with
generator denoted by $x_\al$  
and the stabilizer  $H_\al$,  
and, for $j\geq i$ define a morphism
$X_j \lra X_i$  by sending  $x_j$  to  $g_{ij}x_i$.
These clearly give a surjective $I$-inverse system of non-empty
transitive $G$-sets. By the assumption its inverse limit is non-empty.
 
Suppose $y=(c_ix_i)_{i\in I}\in\invlim X_i$, for some $c_i\in G$.  
By  definitions of $\invlim$ and the maps $X_j\lra X_i$ we have
$g_{ij}(c_jx_i)=c_ix_i$ and, via the isotropy groups $H_i$, this
translates into

$(b)\hfill    g_{ij} + c_j - c_i \in H_i,\,\,i \leq j.\hfill$ 

Note  again  that all the generators of $G$ occurring in (a) have 
both subscripts $\geq\al$, i.e.   $H_\al$  is  contained  in  the
subgroup of $G$ spanned by the generators with this property.  It
thus follows from (b) that $c_i$ and $c_j$ may differ (in the
expansion of $G$) only in the
terms  with  both subscripts in the set $\{i, j, k\}$ of the 
respective generators of $G$, 
$\geq i$.  Hence for every $i<j$ and every $g_{ij}$, 
there is a $k_0>i$ (say $k_0=j$), such that, $\forall k\geq k_0$,
all $c_k$ contain in their expansion the same coefficient
of $g_{ij}$;
in particular, this coefficient will be in all $c_k$ with $k > j$.

If, contrary to our claim, $I$ were of uncountable cofinality, then
this coefficient cannot be nonzero for infinitely (countably) many  
$g_{ij}$, since we could find some $k$
such that $c_k$  involves infinitely many summands in the direct
sum representation, which is
impossible.  Hence only finitely many different $g_{ij}$  have 
nonzero eventual
coefficient, hence we can form an element  $c$  of  $G$  in which each
$g_{ij}$  has this coefficient.  \lqq Translating" our element $y$  of  
$\invlim X_i$
by  $c$,  and redefining the  $c_i$  in terms of this new $y$,  we are
reduced to the situation where the eventual coefficient of
each  $g_{ij}$  is  0.  Hence, by our earlier observations,

$(c)\hfill c_i \text{ involves no }  g_{\al j}  \text{ with }  \al\leq i.\hfill$ 

To complete our proof, let us now map  $G$  homomorphically into the
free abelian group on generators  $\sst{f_i}{i\in I}$ by the 
homomorphism  $D$ defined by

$(d)\hfill  D(g_{ij}) = f_i - f_j.\hfill$ 

Note that $\Ker(D)$  contains all the subgroups  $H_i$.  Hence (b) and (c)
give, respectively:

$(e)\hfill     D(c_i) - D(c_j) = f_i - f_j\,\,\,  (i\leq j),\hfill$ 

$(f)\hfill  D(c_i)  \text{ involves no } f_\al  \text{ with }  \al\leq i.\hfill$ 

\noindent Looking at (e) in the light of (f) we see that the only  $f$  with
subscript  $\leq j$  which  $D(c_i)$  involves is  $f_i$,  and this has
coefficient  1.  Hence fixing  $i$  and taking arbitrarily large  $j$  in
this statement, we conclude that  $f_i$  is the only  $f$  whatsoever that
$D(c_i)$  involves with nonzero coefficient.  Hence

$(g)\hfill  D(c_i) = f_i.\hfill$ 

But from (d) we see that  $D$  carries  $G$  into the subgroup of  $F$  in
which the coefficients of the  $f$'s  sum to  0.  This contradicts (g),
thus the assumption of uncountability of $I$ is false and we conclude
that $I$ either has a maximum element or is of countable cofinality.\qed
\sju

{\bf Notes 1.}
[Bourbaki, 1961,  $\S 3$, Th. 1] attributes to Mittag-Leffler an
implication of the kind (1)$\shimply$(2) where spaces $X_\al$
were taken to be complete metrizable uniform spaces and the
 $I$-inverse systems in (2) 
satisfied the Mittag-Leffler condition, instead of surjective
$I$-inverse systems. 
The proof (3)$\shimply$(1)
and the construction of $E_\al$'s is essentially that of 
[Henkin, 1950] (see also [Bourbaki, 1956, $\S 1$, Exercise 31]).
The result (5)$\shimply$(1) and its proof is by [Bergman, 1998].

\sju

Given an (upward directed) index set $I$, the
{\it surjective cohomological dimension} of $I$ ($scd I$) is the
largest natural number $n$ with $\invlim^{(n)}A_i\neq 0$, for some
surjective $I$-inverse system $\sst{A_i}{i\in I}$ of
sets (or $R$-modules); if no such an $n$ exists it is pronounced to be
$\infty$. Here $\invlim^{(0)}A_i$ is thought of as $\invlim A_i$.
Since we assume to be working within categories with (infinite)
direct products (that can be seen as inverse limits), we get then
that $scd I\geq 0$, for all non-zero $I$.

\proclaim{\quad Theorem 2}
For a non-empty upward directed set $I$, any of the equivalent
statements (1)--(5) in Theorem 1 implies every of the following 
statements:
\roster
\item "(6)" For  every surjective $I$-inverse system of non-trivial free 
abelian  groups  (modules),  its  inverse   limit   is   likewise
non-trivial.

\item "(7)" For  every  surjective  $I$-inverse  system  $\sst{M_\al,
f_{\al\be}}{\al,\be\in I}$ of  non-trivial 
abelian  groups  (modules),  its  inverse   limit   is   likewise
non-trivial. 

\item "(8)" For  every  surjective  I-inverse  system $\bold A$ of 
abelian groups (modules) and every $n\geq 1$, $\invlim^{(n)}A_i=0$. 

\item "(9)" For  every  surjective  I-inverse  system $\bold A$ of
abelian groups (modules) and 
every exact sequence $0\lra \bold A\lra \bold B\lra \bold C\lra 0$ 
of $I$-inverse systems,  the
corresponding sequence of inverse limits 
$0\lra \invlim\bold A\lra \invlim\bold B\lra \invlim\bold C\lra 0$ 
is likewise exact.

\item "(10)" $scd I=0$. 
\endroster
\endproclaim

{\sl Proof.} Note first that (8)$\iff$(9)$\iff$(10), since $scd
I\geq 0$. For (8)$\iff$(9)
use the exact sequence $0\lra\invlim A_i\lra\invlim
B_i\lra\invlim C_i\lra\invlim^{(1)} A_i\lra\dots $

(7)$\shimply$(6): since the former is a more general
than the latter. 

(1)$\shimply$(7): The proof of this is same, mutatis mutandis, as
(1)$\shimply$(2); we pick the starting element $x_{\al_0}\in
M_{\al_0}$  (or $x_{\al_1}\in M_{\al_1}$) to be non-zero, to
ensure that the resulting element in the limit is also non-zero.

(2)$\shimply$(9): Let $\bold A=\{A_j,f_{ij}\}$, $\bold
B=\{B_j,g_{ij}\}$, $\bold C=\{C_j,h_{ij}\}$, and exact sequences 

$$(*)\qquad
\CD
0\lra A_j @>{u_j} >> B_j @>{v_j} >>C_j\lra 0 \\
\qquad\quad @V{f_{ij}}VV   @V{g_{ij}}VV   \hskip-1truecm@Vh_{ij}VV \\
0\lra A_i @>{u_i} >> B_i @>{v_i} >>C_i\lra 0 
\endCD
$$

\noindent Denote also $v=\invlim v_i$, $u=\invlim u_i$. Given
$0\neq c=(c_i)_{i\in I}\in\invlim C_i$ (thus $h_{ij}(c_j)=c_i$)
we need a $b=(b_i)_{i\in I}\in\invlim B_i$ with $v(b)=c$, i.e.
$\forall j\,\, v_j(b_j)=c_j$. Denote $E_j=v_i^{-1}(c_j)\neq\ets$;
it is non-empty, since the $v_j$'s are surjective. Denote now
$e_{ij}=g_{ij}|E_j$; it is straightforward to show that 
$\bold E=\{E_j, e_{ij}\}$ is an $I$-inverse system. We now show
that it is surjective: To this end, start with a $b_i\pr\in E_i$,
i.e. such that $v_i(b\pr_i)=0$. By exactness in (*), there is an
$a_i\in A_i$ with $u_i(a_i)=b\pr_i$ (**). By surjectivity of the $f$'s,
there is an $a_j\in A_j$ with $f_{ij}(a_j)=a_i$ (***). 
Since $u_j(a_j)\in\ima u_j=\Ker v_j$, we have $v_j(u_j(a_j))=0$, i.e.
$u_j(a_j)\in E_j$. Appeal again to (*), then (***) and (**) to get 
$g_{ij}u_j(a_j)=u_if_{ij}(a_j)=u_i(a_i)=b_i\pr$; this proves
surjectivity of all $e_{ij}$. By the assumption, $\invlim
E_j\neq\ets$, hence for any $b$ in that set we have $v(b)=c$ by
the very construction, which proves the claim.  \qed

{\bf Notes 2.} By way of universal images, an aditional set of (equivalent)
statements may be added with the word \lqq surjective" replaced by the
word the ML condition, in Theorem 2.
(1)$\shimply$(9) was proved in [Grothendieck, 1961] 
where $\bold A$ is required to satisfy the
Mittag-Leffler condition; 
[Goblot, 1970] replaces countability of the index set $I$ by the
requirement that $I$ is well-ordered and that the participating
objects and maps form a continuous (smooth) system. [Jensen, 1972] replaces
countability of $I$ by a requirement that $\bold A$ is semi-flabby.
In an ongoing work we will show that $scd I=n$ if and only if
$|I|=\aleph_n$, otherwise the surjective cohomological dimension 
is infinite ([Mitchell, 1973] shows similar result for the
cohomological dimension: For a directed $I$, if $cf I=\aleph_k$,
$-1\leq k\leq\infty$, then $cd I=n+1$ iff $k=n$).  This will
then establish equivalence of conditions in Theorem 1 with
conditions in Theorem 2. 

\ju

\centerline{\bf References}
\sju
\noindent [1] George Bergman, Private communication, Mar 5, 1998, originated
likely in 1992

\noindent [2] N. Bourbaki, {\it Topologie g\'en\'erale,} chap. II, 3$^e$ \'ed,  1961

\noindent [3] N. Bourbaki, {\it Th\'eorie des ensembles}, chap. III, 1956


\noindent [4] Goblot, R\'emi
Sur les d\'eriv\'es de certaines limites projectives.
Applications aux modules.
Bull. Sc. math., France, 2nd series, 94(1970), 251--255

\noindent [5] A. Grothendieck, \'El\'ements de g\'eom\'etrie alg\'ebrique, III,
Pr\'eliminaires
{\it Inst. Hautes \'Etudes Sci. Publ. Math.}, {\bf 11}(1961), 343--423 (1--79)

\noindent [6] Leon Henkin, A problem on inverse mapping systems, {\it Proc.
Amer. Math. Soc.} {\bf 1}(1950), 224-225

\noindent [7] C.U. Jensen, Les Foncteurs D\'eriv\'es de $\invlim$ et leurs
Applications en Th\'eorie des Modules, Lecture Notes in
Mathematics,  Springer-Verlag, Berlin, 1972

\noindent [8] Barry Mitchell, The cohomological dimension of a directed set,
{\it Can. J. Math. } {\bf 25}(1973), No.2, 233--238
\vfill

\enddocument